\documentclass[10pt]{amsart}

\usepackage[all, cmtip]{xy}
\usepackage{amsmath,amsthm,amssymb,amsfonts,amscd}
\usepackage{pstricks}
 \usepackage{epsfig}
 \usepackage{pst-grad} % For gradients
 \usepackage{pst-plot} % For axes

\usepackage{graphicx}
\newtheorem{theorem}{Theorem}[section]
\newtheorem{lemma}[theorem]{Lemma}
\newtheorem{proposition}[theorem]{Proposition}

\theoremstyle{definition}

\newtheorem{notation}[theorem]{Notation}

\newtheorem{conjecture}[theorem]{Conjecture}

\theoremstyle{remark}
\newtheorem{remark}[theorem]{Remark}

\numberwithin{equation}{section}

\begin{document}

\title[rational surfaces with ample canonical divisor]{Constructions of singular rational surfaces of Picard number one with ample canonical divisor}

\author{DongSeon Hwang}
\address{School of Mathematics, Korea Institute for Advanced Study, Seoul 130-722, Korea}
%\curraddr{Department of Mathematics and Statistics,
%Case Western Reserve University, Cleveland, Ohio 43403}
\email{dshwang@kias.re.kr}
\thanks{Research supported by the National Research Foundation of Korea (NRF-2007-C00002).}

%\author{Author Two}
\author{JongHae Keum}
%\address{School of Mathematics, Korea Institute For Advanced Study, Seoul 130-722, Korea}
\email{jhkeum@kias.re.kr}

%    General info
%\subjclass[2010]{Primary 14J17, 14E20; Secondary 46E25, 20C20}

\date{July 8, 2010}
%\date{June 22, 2009 and, in revised form,}



\begin{abstract} Koll\'ar gave a series of examples of rational surfaces of Picard number $1$ with ample canonical divisor having cyclic
singularities. In this paper, we construct several series of new
examples in a geometric way, i.e., by blowing up several times
inside a configuration of curves on the projective plane and then
by contracting chains of rational curves. One series of our
examples have the same singularities as Koll\'ar's examples.
\end{abstract}

\maketitle

%\section*{This is an unnumbered first-level section head}
%This is an example of an unnumbered first-level heading.

%% The correct journal style for \specialsection is all uppercase; a known bug
%% in amsart.cls prevents this, so input must be uppercase until it is fixed.
%\specialsection*{This is a Special Section Head}
%\specialsection*{THIS IS A SPECIAL SECTION HEAD}
%This is an example of a special section head%
%%%%%%%%%%%%%%%%%%%%%%%%%%%%%%%%%%%%%%%%%%%%%%%%%%%%%%%%%%%%%%%%%%%%%%%%
%\footnote{Here is an example of a footnote. Notice that this footnote
%text is running on so that it can stand as an example of how a footnote
%with separate paragraphs should be written.
%\par
%And here is the beginning of the second paragraph.}%
%%%%%%%%%%%%%%%%%%%%%%%%%%%%%%%%%%%%%%%%%%%%%%%%%%%%%%%%%%%%%%%%%%%%%%%%

\section{Introduction}

A rational surface $S$ with quotient singularities has been
studied extensively when its anti-canonical divisor $-K_S$ is
ample or numerically trivial. In the former case the surface is
called a log del Pezzo surface, and in the latter case the surface
is called a log Enriques surface. On the other hand, when $K_S$ is
ample, very little is known about the classification of such
surfaces. Moreover, if in addition $S$ has Picard number
$\rho(S)=1$, nothing seems to be known except the examples due to
Koll\'ar (\cite{Kol08}, Example 43).

Koll\'ar constructed a series of such examples by contracting two
rational curves on some well-chosen weighted projective
hypersurfaces. We briefly review his construction. Let
$$Y = Y(a_1, a_2, a_3, a_4) := (x_1^{a_1}x_2 + x_2^{a_2}x_3 + x_3^{a_3}x_4 + x_4^{a_4}x_1 = 0)$$
be a hypersurface in $\mathbb{P}(w_1, w_2, w_3, w_4)$, where $a_i$
and the weights $w_i$ satisfy a system of equations
$$a_1w_1+w_2 = a_2w_2 + w_3 = a_3w_3 + w_4 = a_4w_4+w_1 = d$$
with solutions
$$\begin{array}{ll} w_1 = \frac{1}{w^*}(a_2 a_3 a_4-a_3 a_4+a_4-1),
& w_2 = \frac{1}{w^*}(a_1 a_3 a_4-a_1 a_4+a_1-1),\\
  w_3 =
\frac{1}{w^*}(a_1 a_2 a_4-a_1 a_2 + a_2-1), &  w_4 =
\frac{1}{w^*}(a_1 a_2 a_3- a_2 a_3 + a_3-1),
\end{array}$$
$$d = \frac{1}{w^*}(a_1 a_2 a_3 a_4-1)$$ where $w^* = \gcd(w_1, w_2,
w_3, w_4).$ If $w^* = 1$, then by Theorem 39 in \cite{Kol08}, $Y$
is a rational surface with $4$ cyclic singularities at the
coordinate vertices and with $H_2(Y, \mathbb{Q}) \cong
\mathbb{Q}^3$. The two rational curves
$$C_1 := (x_1 = x_3 = 0),\quad C_2 := (x_2 = x_4 = 0)$$
 are extremal rays
for the $K_Y + (1- \epsilon)(C_1 + C_2)$ minimal model program for
$0 < \epsilon \ll 1$. Thus $C_1$ and $C_2$  are both contractible
to quotient singularities and we get a rational surface of Picard
number 1,
$$\pi : Y=Y(a_1,a_2,a_3,a_4) \rightarrow X = X(a_1, a_2, a_3, a_4).$$
If $a_1, a_2, a_3, a_4 \geq 4$, then $K_{X}$ is ample by Theorem
39(5) in \cite{Kol08}. The surface $X$ has two cyclic
singularities and no other singularities.

First we determine the types of singularities of $X$.

\begin{theorem}\label{X} Assume that $w^* = 1$.
Then the two cyclic singularities of the surface $X(a_1, a_2, a_3,
a_4)$ are of type
$$\frac{1}{s_1}(w_2,w_4)\,\,\,{\rm and}\,\,\,\, \frac{1}{s_2}(w_1,w_3)$$
where $s_1=a_4w_4-w_3$ and $s_2=a_3w_3-w_2$. Their Hirzebruch-Jung
continued fractions are $$[\underset{a_4-1}{\underbrace{2, \ldots,
2}}, a_3, a_1, \underset{a_2-1}{\underbrace{2, \ldots,
            2}}]\,\,\,{\rm and}\,\,\,\, [\underset{a_3-1}{\underbrace{2, \ldots, 2}},
            a_2, a_4, \underset{a_1-1}{\underbrace{2, \ldots,
            2}}],$$
            respectively.
\end{theorem}

\begin{remark} The condition $w^* = 1$ gives some restriction on
the choice of the integers $a_i$ for $X(a_1, a_2, a_3, a_4)$. In particular, it rules out the possibility that $a_1=a_2=a_3=a_4$.
\end{remark}

We recall that a normal projective surface $S$ with the same Betti
numbers with $\mathbb{P}^2$ is called a
\emph{$\mathbb{Q}$-homology projective plane}. When a normal
projective surface $S$ has quotient singularities only, $S$ is a
$\mathbb{Q}$-homology projective plane if the second Betti number
$b_2(S)=1$. For convenience, we adopt the following terminology: a
$\mathbb{Q}$-homology projective plane which is a rational surface
is called a \emph{rational $\mathbb{Q}$-homology projective
plane}. The surfaces $X(a_1, a_2, a_3, a_4)$ are rational
$\mathbb{Q}$-homology projective planes if $w^* = 1$.

Next, we  give a geometric construction of a series of rational
$\mathbb{Q}$-homology projective planes with ample canonical
divisor having the same singularities as Koll\'ar's examples.

\begin{theorem}\label{two} For each integers $a_1, a_2, a_3, a_4\ge 2$,
there exists a rational $\mathbb{Q}$-homology projective plane
$T=T(a_1, a_2, a_3, a_4)$ with two cyclic singularities of type
$$[\underset{a_4-1}{\underbrace{2, \ldots,
2}}, a_3, a_1, \underset{a_2-1}{\underbrace{2, \ldots,
            2}}]\,\,\,{\rm and}\,\,\,\, [\underset{a_3-1}{\underbrace{2, \ldots, 2}},
            a_2, a_4, \underset{a_1-1}{\underbrace{2, \ldots,
            2}}].$$
The surfaces can be constructed by blowing up several times inside
general $4$ lines of $\mathbb{P}^2$ and then by contracting two
chains of rational curves. Moreover,
    \begin{enumerate}
    \item if $a_1, a_2, a_3, a_4\ge 3$ and $a_i >3$ for some
        $i$, then $K_T$ is ample,
    \item if $a_1=a_2=a_3=a_4= 3$, then $K_T$ is
        numerically trivial.
    \end{enumerate}
\end{theorem}

\begin{remark} Note that for any choice of the numbers $a_i\ge 2$, $T(a_1, a_2, a_3, a_4)$ exists.
Thus in terms of types of singularities our examples $T(a_1, a_2,
a_3, a_4)$ include properly Koll\'ar's examples $X(a_1, a_2, a_3,
a_4)$.
\end{remark}

 Starting with different configurations of curves
on $\mathbb{P}^2$ (see Section 5 and 6 for construction), we also
construct other series of rational $\mathbb{Q}$-homology
projective planes with $K_S$ ample having one or three cyclic
singularities.

\begin{theorem}\label{one}  For each integer $b \geq 2$,
there exists a rational $\mathbb{Q}$-homology projective plane
$S:=S(b)$ with a unique cyclic singularity of type
$$\frac{1}{27b^2-36b+4}(1, 9b^2-9b+1).$$ Moreover,
    \begin{enumerate}
    \item if $b=2$, then $K_S$ is numerically trivial,
    \item if $b>2$, then $K_S$ is ample.
    \end{enumerate}
\end{theorem}

\begin{theorem}\label{three}  For each integer $b \geq 2$,
there exists a rational $\mathbb{Q}$-homology projective plane
$S:=S(b)$ with three cyclic singularities of type
$$A_1, \,\,\, \frac{1}{7}(1,3),\,\,\, \frac{1}{3b^2-2b-2}(1, 2b^2-b-1).$$
Moreover,
    \begin{enumerate}
    \item if $b<5$, then $-K_S$ is ample,
    \item if $b=5$, then $K_S$ is numerically trivial,
    \item if $b>5$, then $K_S$ is ample.
    \end{enumerate}
\end{theorem}

Using different configurations of curves on $\mathbb{P}^2$ or
different blow-ups of the same configuration, it is possible to
construct more examples with $K_S$ ample having at most three
cyclic singularities.

  The authors have shown that every
rational $\mathbb{Q}$-homology projective plane with quotient
singularities has at most $4$ singularities \cite{HK1}. It seems
impossible, or very difficult, to construct a rational
$\mathbb{Q}$-homology projective plane with $K_S$ ample having 4
cyclic singularities. We recall that there are rational
$\mathbb{Q}$-homology projective planes with $K_S$ ample having an
arbitrary number of rational singularities (\cite{HK1},
Introduction).

%, and there are rational $\mathbb{Q}$-homology projective planes with $K_S$ ample having 4 quotient singularities, not all cyclic.

%\begin{conjecture}\label{amy}
%Let $S$ be a rational $\mathbb{Q}$-homology projective plane with cyclic singularities. If $K_S$ is ample then $S$ has at most $3$ singular points. \end{conjecture}

The original motivation of this study is the following conjecture called
the algebraic Montgomery-Yang problem:

\begin{conjecture}[\cite{Kol08}, Conjecture 30]\label{amy}
 Let $S$ be a
$\mathbb{Q}$-homology projective plane with quotient
singularities. If $\pi_1(S^0) = 1$, then $S$ has at most $3$
singular points. Here $S^0$ denotes the smooth locus of $S$.
\end{conjecture}

 Recently, the authors have confirmed Conjecture \ref{amy} unless
 $S$ is a rational $\mathbb{Q}$-homology projective plane with cyclic
singularities having ample canonical divisor (\cite{HK2},
\cite{HK3}). Therefore, to solve the conjecture we need to study
such surfaces. The result in this paper is the first step toward
the goal.

\medskip
 Throughout this paper, we work over the field $\mathbb{C}$ of complex numbers.

\bigskip
{\bf Acknowledgements.} DongSeon Hwang thanks Stephen Coughlan for
helpful conversations. Both authors thank Miles Reid for helpful
comments.

\bigskip

\section{Preliminaries}

 Let $\mathcal{H}$ be the set of all Hirzebruch-Jung continued
fractions
 $$[n_1, n_2, ..., n_l]= n_1 - \dfrac{1}{n_2-\dfrac{1}{\ddots -
\dfrac{1}{n_l}}},$$ i.e.,
$$\mathcal{H} = \underset{l \geq 1}{\bigcup} \{ [n_1, n_2, \ldots, n_l] \mid \textrm{all}\,\,
n_j\,\, \textrm{are integers} \geq 2 \}.$$   We collect some
notations and properties of Hirzebruch-Jung continued fractions
for later use.

\begin{notation}\label{uv-n}
    For a fixed $w = [n_1, n_2, \ldots, n_l] \in \mathcal{H}$, we define
        \begin{enumerate}
        \item $|w|=q$, the order of
the cyclic singularity corresponding to $w$, i.e.,
$w=\dfrac{q}{q_1}$ with $1\le q_1<q, \,\, \gcd(q, q_1)=1$. Note
that $|w|$ is the absolute value of the determinant of the
intersection matrix corresponding to $w$. We also define
\item   $u_j :=  |[n_1, n_2, \ldots, n_{j-1}]|\,\,\,
                (2\leq j\leq l+1), \quad u_0=0, \,\,u_1=1$,
 \item   $v_j := |[n_{j+1}, n_{j+2}, \ldots, n_l]|
                \,\,\, (0\leq j\leq l-1), \quad v_l=1,
                \,\,v_{l+1}=0$.
 \end{enumerate}
        Note that   $u_{l+1}=v_0=|w|$.
\end{notation}

\begin{lemma}\label{hk}$($Lemma 2.4 (6) in \cite{HK3}$)$
 $$|[n_1,\ldots, n_{j-1}, n_j+1, n_{j+1}, \ldots,
        n_l]|= v_j u_j +|[n_1, n_2, \ldots, n_l]|.$$
\end{lemma}

\begin{notation}
    $[2*a,b,c,2*d] := [\underset{a}{\underbrace{2, \ldots, 2}},
            b, c, \underset{d}{\underbrace{2, \ldots,
            2}}]$
\end{notation}

\begin{lemma}\label{cf}
$$|[2*(a-1),b,c,2*(d-1)]| =
        abcd-abd-acd+ab+cd-a-d+1.$$
\end{lemma}

\begin{proof}
It is easy to see that $$|[2*(a-1), 2, 2, 2*(d-1)]| = a+d+1.$$ By
applying Lemma \ref{hk} $(b-2)$ times, we get
$$|[2*(a-1), b, 2, 2*(d-1)]| = (a+d+1) +
(d+1)a(b-2).$$ Again, by applying Lemma \ref{hk} $(c-2)$ times, we
get $$|[2*(a-1), b, c, 2*(d-1)]| = (a+d+1) + (d+1)a(b-2) +
d(ab-a+1)(c-2),$$ from which the result follows.
\end{proof}

Now let $X$ be a $\mathbb{Q}$-homology projective plane with only
cyclic singularities, and $f:X' \rightarrow X$ be its minimal
resolution. For each singular point $p$ of $X$, the dual graph of
$f^{-1}(p)$ is of the form
$$\underset{A_{1,p}}{\overset{-n_{1,p}}\circ}
-\underset{A_{2,p}}{\overset{-n_{2,p}}\circ}-\cdots
-\underset{A_{l_p,p}}{\overset{-n_{l_p,p}}\circ}$$ which
corresponds to the Hirzebruch-Jung continued fraction
$$w_p = [n_{1,p}, n_{2, p}, \ldots, n_{l_p, p}]$$
where $A_{j, p}$'s are irreducible components of   $f^{-1}(p)$. We
omit the subscript $p$ if the meaning is clear from the context.

\begin{lemma}\label{ample}
Let $E$ be a $(-1)$-curve on $X'$. Then
$$  Ef^*(K_X) =-1+  \underset{p \in Sing(S)}{\sum}
\overset{l_p}{\underset{j = 1}{\sum}} \big( 1 - \dfrac{v_{j,p} +
      u_{j,p}}{|w_p|} \big) EA_{j,p}.$$
      \end{lemma}

\begin{proof}
By Lemma 2.2 of \cite{HK1}, we have the adjunction formula
$$K_{X'}=f^*(K_X)-\underset{p \in Sing(S)}{\sum}
\overset{l_p}{\underset{j = 1}{\sum}} \big( 1 - \dfrac{v_{j,p} +
      u_{j,p}}{|w_p|} \big) A_{j,p}$$
Intersecting the adjunction formula with $E$ we get the equality.
\end{proof}

We have the following criterion for ampleness of $K_X$.

\begin{lemma}\label{ample2}
      \begin{enumerate}
      \item $K_X$ is ample iff $Ef^*(K_X)>0$ for all irreducible curves $E$ not contracted by $f$
      iff $Ef^*(K_X)>0$ for an irreducible curve $E$ not contracted by $f$.
      \item $K_X$ is numerically trivial iff $Ef^*(K_X)=0$ for all irreducible curves $E$ not contracted by $f$
      iff $Ef^*(K_X)=0$ for an irreducible curve $E$ not contracted by $f$.
      \item $-K_X$ is ample iff $Ef^*(K_X)<0$ for all irreducible curves $E$ not contracted by $f$
      iff $Ef^*(K_X)<0$ for an irreducible curve $E$ not contracted by $f$.
      \end{enumerate}
\end{lemma}

\begin{proof}
Since $X$ has Picard number 1, we have the trichotomy: $K_X$ is
ample; $K_X$ is numerically trivial; $-K_X$ is ample. The
assertion follows from $Ef^*(K_X)=f_*(E)K_X$.
\end{proof}

\section{Singularities of Koll\'ar's examples}

In this section, we will prove Theorem \ref{X} by using the
technique of unprojection \cite{Reid}. Let
$$F = x_1^{a_1}x_2 + x_2^{a_2}x_3 + x_3^{a_3}x_4 + x_4^{a_4}x_1$$
 be the weighted homogeneous polynomial defining $$Y=Y(a_1,a_2,a_3,a_4)\subset\mathbb{P}(w_1, w_2, w_3, w_4).$$
%\subsection{Construction via unprojection} In this subsection we interpret \\Koll\`ar's examples using unprojection. We refer the reader to \cite{Reid} for the the theory of unprojection; Section 2 in \cite{NP1} and \cite{NP2} for the parallel unprojection.
Recall that on $Y$ there are two disjoint rational curves $C_1 =
\mathbb{P}(w_2, w_4)$ and $C_2 = \mathbb{P}(w_1, w_3)$. Write $F$
as
$$F = A x_1 + B x_3$$ where $$A = x_1^{a_1-1} x_2 + x_4^{a_4}, \quad B = x_2^{a_2} + x_3^{a_3-1}x_4.$$ By introducing a new variable
$$y_1 = \frac{A}{x_3} = - \frac{B}{x_1},$$
we get an unprojection morphism $Y \rightarrow Y_1$ where
$$Y_1  =   (x_3 y_1 = A, x_1 y_1 = -B) \subset \mathbb{P}(w_1, w_2, w_3, w_4, s_1),$$
where $$s_1 := \deg(y_1)= a_4w_4-w_3=a_2w_2-w_1.$$
 This morphism contracts $C_1$ to the singular point
 $$p=(0,0,0,0,1)\in Y_1$$ of $Y_1$.
It is easy to see that near $p$ $$Y_1\cong \mathbb{P}(w_2, w_4,
s_1).$$ Thus the singular point $p$ is of type $\frac{1}{s_1}(w_2,
w_4)$.

Similarly, using another unprojection morphism contracting $C_2$,
we see that the other singularity is of type $\frac{1}{s_2}(w_1,
w_3)$ where $$s_2= a_1w_1-w_4=a_3w_3-w_2.$$

%Write $F$ by $F = A' x_2 + B' x_4$ where $A' =x_1^{a_1} + x_2^{a_2-1} x_3$ and $B' = x_3^{a_3} + x_4^{a_4-1}x_1.$ By introducing a new unprojection variable
%$$t_2 = \frac{A'}{x_4} =-\frac{B'}{x_2}$$ with $s_2 := deg (t_2)$, we get the unprojection morphism $S \rightarrow S_2$  where
%$$S_2  =   (x_4 t_2 = A', x_2 t_2 = -B') \subset \mathbb{P}(w_1, w_2, w_3, w_4, s_2),$$ which contracts $C_2$ to a singular point of type $\frac{1}{s_2}(w_1, w_3)$.

Now we determine the types of singularities in terms of
Hirzebruch-Jung continued fraction.
\begin{notation}
The expression $\gamma\underset{n}{\equiv}\delta$ means that
$\gamma$ is congruent to $\delta$ modulo $n$.
\end{notation}

Write
$$\frac{1}{s_1}(w_2, w_4) = \frac{1}{s_1}(1, t_1), \qquad
\frac{1}{s_2}(w_1, w_3) = \frac{1}{s_2}(1, t_2)$$ where $t_1
\underset{s_1}{\equiv}  \alpha w_4$ and $t_2
\underset{s_2}{\equiv} \beta w_3$ for some integers $\alpha$ and
$\beta$ satisfying
$$\alpha w_2 \underset{s_1}{\equiv}  1, \qquad \beta w_1 \underset{s_2}{\equiv}  1.$$

The following lemma is immediate.

\begin{lemma}\label{mod}
For an integer $t$,
    \begin{enumerate}
    \item $t \underset{s_1}{\equiv} \alpha w_4$ if and only if $tw_2 \underset{s_1}{\equiv} w_4,$
    \item $t \underset{s_2}{\equiv} \beta w_3$ if and only if $tw_1 \underset{s_2}{\equiv} w_3.$
    \end{enumerate}
\end{lemma}

\begin{lemma}\label{type} If $w^* = 1$,  then the following equalities hold:
    \begin{enumerate}
        \item $[2*(a_4-1),
            a_3, a_1, 2*(a_2-1)] = \dfrac{s_1}{t_1}$.
        \item $[2*(a_3-1),
            a_2, a_4, 2*(a_1-1)] = \dfrac{s_2}{t_2}$.
    \end{enumerate}
\end{lemma}

\begin{proof}
 (1) Note that
$$\begin{array}{lll}
s_1 &=& a_4w_4-w_3\\ &=& a_1a_2a_3a_4 - a_1a_2a_4 - a_2a_3a_4 +a_1a_2
+ a_3a_4 -a_2-a_4+1\\
& =&|[2*(a_4-1), a_3, a_1, 2*(a_2-1)]|
\end{array}$$
where the last equality follows from Lemma \ref{cf}.\\ Now it is
enough to show that $|[2*(a_4-2), a_3, a_1, 2*(a_2-1)]| = t_1$,
i.e.,
$$|[2*(a_4-2), a_3, a_1,
2*(a_2-1)]|  \underset{s_1}{\equiv} \alpha w_4.$$ Since
$$w_2 = a_1 a_3 a_4-a_1 a_4+a_1-1\,\,\,{\rm and}\,\,\, w_4 =a_1 a_2 a_3- a_2 a_3 +
a_3-1,$$ a direct computation shows
$$\begin{array}{lll} |[2*(a_4-2), a_3, a_1, 2*(a_2-1)]| w_2 &=& w_4
+ (a_1a_3a_4-a_1a_3-a_1a_4+2a_1-1)s_1\\&\underset{s_1}{\equiv}&
w_4.\end{array}$$ Now the assertion follows from Lemma \ref{mod}.

(2) Similarly, note that\\
$$\begin{array}{lll} s_2 &=& a_1w_1-w_4\\ &=& a_1a_2a_3a_4 -
a_1a_2a_3 - a_1a_3a_4 +a_1a_4
+ a_2a_3 -a_1-a_3+1\\
& =&|[2*(a_3-1), a_2, a_4, 2*(a_1-1)]|.
\end{array}$$
A direct computation also shows that
$$\begin{array}{lll} |[2*(a_3-2), a_2, a_4, 2*(a_1-1)]| w_1 &=& w_3 +
(a_2a_3a_4-a_2a_4-a_3a_4+2a_4-1)s_2\\
&\underset{s_1}{\equiv}& w_3.\end{array}$$ Again the assertion
follows from Lemma \ref{mod}.
\end{proof}

We have proved Theorem \ref{X}.

\section{Construction of $T(a_1, a_2, a_3, a_4)$}

We construct rational $\mathbb{Q}$-homology projective planes with
ample canonical divisor having the same singularities as
Koll\`ar's examples starting from the configuration of general
four lines on $\mathbb{P}^2$.

Consider four general lines $L_1, L_2, L_3, L_4$ on
$\mathbb{P}^2$. Choose four points among the six intersection
points such that each $L_i$ passes through two of them.

\bigskip

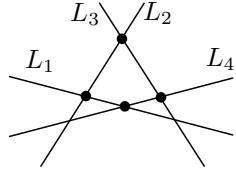
\scalebox{1} % Change this value to rescale the drawing.
{
\begin{pspicture}(-4, -1)(5,1.5)
\psline[linewidth=0.02cm](1.2185937,1.0816406)(2.6185937,-1.0983593)
\psline[linewidth=0.02cm](1.8185937,1.0816406)(0.43859375,-1.0983593)
\psline[linewidth=0.02cm](0.01859375,-0.6983594)(2.9985938,0.0816406)
\psline[linewidth=0.02cm](0.01859375,0.1016406)(2.9985938,-0.6783594)
\psdots[dotsize=0.14](1.5585938,-0.29835936)
\psdots[dotsize=0.14](2.0385938,-0.17835934)
\psdots[dotsize=0.14](1.0385938,-0.15835935)
\psdots[dotsize=0.14](1.5185938,0.6016406)
\usefont{T1}{ptm}{m}{n}
\rput(0.43234375,0.30992186){$L_1$}
\usefont{T1}{ptm}{m}{n}
\rput(1.0323437,0.92992187){$L_3$}
\usefont{T1}{ptm}{m}{n}
\rput(2.0,0.9499219){$L_2$}
\usefont{T1}{ptm}{m}{n}
\rput(2.8323438,0.32992184){$L_4$}
\end{pspicture}
}

\bigskip\noindent

By blowing up each of the four marked intersection points twice,
we get a rational surface $Z(2,2,2,2)$ having 12 rational curves
such that

\bigskip

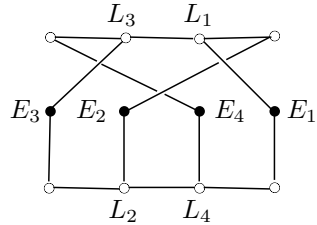
\scalebox{1} % Change this value to rescale the drawing.
{
\begin{pspicture}(-3.5, -1.5)(5,1.5)
\psdots[dotsize=0.14,fillstyle=solid,dotstyle=o](0.7650486,1.0255765)
\psdots[dotsize=0.14,fillstyle=solid,dotstyle=o](1.7650486,1.0204202)
\psdots[dotsize=0.14,fillstyle=solid,dotstyle=o](2.7450485,1.0055765)
\psdots[dotsize=0.14,fillstyle=solid,dotstyle=o](3.7450485,1.0455763)
\psdots[dotsize=0.14,fillstyle=solid,dotstyle=o](0.7450486,-0.9795798)
\psdots[dotsize=0.14,fillstyle=solid,dotstyle=o](1.7450486,-0.9795797)
\psdots[dotsize=0.14,fillstyle=solid,dotstyle=o](2.7450485,-0.9795797)
\psdots[dotsize=0.14,fillstyle=solid,dotstyle=o](3.7450485,-0.9744236)
\psdots[dotsize=0.14](0.7650486,0.040420208)
\psdots[dotsize=0.14](1.7450486,0.040420208)
\psdots[dotsize=0.14](2.7450485,0.040420208)
\psdots[dotsize=0.14](3.7450485,0.040420208)
\psline[linewidth=0.02cm](1.7409375,0.0296875)(1.7409375,-0.8903125)
\psline[linewidth=0.02cm](2.7409377,0.0496875)(2.7409377,-0.8903125)
\psline[linewidth=0.02cm](3.7409377,0.0496875)(3.7409377,-0.8703125)
\psline[linewidth=0.02cm](0.7609375,0.0496875)(0.7409375,-0.8903125)
\psline[linewidth=0.02cm](2.8009374,0.96968746)(3.7409377,0.0296875)
\psline[linewidth=0.02cm](0.7409375,0.0496875)(1.7209375,0.96968746)
\psline[linewidth=0.02cm](1.7409375,0.0296875)(3.0009375,0.6896875)
\psline[linewidth=0.02cm](3.1009374,0.7296875)(3.6809375,1.0296875)
\psline[linewidth=0.02cm](0.8209375,0.9896875)(1.3809375,0.7096875)
\psline[linewidth=0.02cm](1.4409375,0.6696875)(2.1809375,0.3096875)
\psline[linewidth=0.02cm](2.2809374,0.2696875)(2.7209377,0.0496875)
\usefont{T1}{ptm}{m}{n}
\rput(2.7323437,1.3196875){$L_1$}
\usefont{T1}{ptm}{m}{n}
\rput(1.7123437,1.3196875){$L_3$}
\usefont{T1}{ptm}{m}{n}
\rput(1.7323438,-1.2803125){$L_2$}
\usefont{T1}{ptm}{m}{n}
\rput(2.7123437,-1.2803125){$L_4$}
\usefont{T1}{ptm}{m}{n}
\rput(4.112344,0.040420208){$E_1$}
\usefont{T1}{ptm}{m}{n}
\rput(1.3123437,0.040420208){$E_2$}
\usefont{T1}{ptm}{m}{n}
\rput(0.43234375,0.040420208){$E_3$}
\usefont{T1}{ptm}{m}{n}
\rput(3.1523438,0.040420208){$E_4$}
\psline[linewidth=0.02cm](0.8209375,1.0296875)(1.7009375,1.0096875)
\psline[linewidth=0.02cm](1.8209375,1.0296875)(2.6609375,1.0096875)
\psline[linewidth=0.02cm](2.8009374,1.0096875)(3.7009375,1.0296875)
\psline[linewidth=0.02cm](0.8009375,-0.9703125)(1.6809375,-0.9903125)
\psline[linewidth=0.02cm](1.8009375,-0.9703125)(2.6809375,-0.9703125)
\psline[linewidth=0.02cm](2.8009374,-0.9703125)(3.6809375,-0.9903125)
\end{pspicture}
}

\bigskip\noindent
is their dual graph. Here, $\bullet$  means a $(-1)$-curve and
$\circ$ means a $(-2)$-curve.

Let $Z(2+r_1, 2+r_2, 2+r_3, 2+r_4)$ be the surface obtained from
$Z(2,2,2,2)$ by blowing up  $r_1$ times at  $E_1 \cap L_1$, $r_2$
times at $E_2 \cap L_2$, $r_3$ times at $E_3 \cap L_3$  and $r_4$
times at $E_4 \cap L_4$, respectively. Set $$(a_1,
a_2,a_3,a_4):=(2+r_1, 2+r_2, 2+r_3, 2+r_4).$$ Then $Z(a_1,
a_2,a_3,a_4)$ has the following configuration of rational curves.

\bigskip
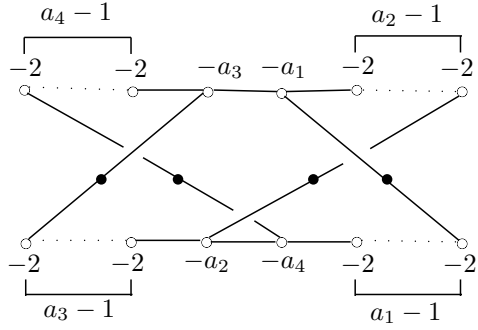
\scalebox{1} % Change this value to rescale the drawing.
{
\begin{pspicture}(-2.5, -2)(5,2)
\psdots[dotsize=0.14,fillstyle=solid,dotstyle=o](1.7450486,0.97557646)
\psdots[dotsize=0.14,fillstyle=solid,dotstyle=o](2.7450485,0.9704202)
\psdots[dotsize=0.14,fillstyle=solid,dotstyle=o](3.7250485,0.9555765)
\psdots[dotsize=0.14,fillstyle=solid,dotstyle=o](4.7250485,0.9955764)
\psdots[dotsize=0.14,fillstyle=solid,dotstyle=o](1.7250487,-1.0295798)
\psdots[dotsize=0.14,fillstyle=solid,dotstyle=o](2.7250485,-1.0295798)
\psdots[dotsize=0.14,fillstyle=solid,dotstyle=o](3.7250485,-1.0295798)
\psdots[dotsize=0.14,fillstyle=solid,dotstyle=o](4.7250485,-1.0244236)
\psdots[dotsize=0.14](1.3250486,-0.18957978)
\psdots[dotsize=0.14](2.3450487,-0.18957978)
\psdots[dotsize=0.14](4.1450486,-0.18957978)
\psdots[dotsize=0.14](5.1250486,-0.18957978)
\psdots[dotsize=0.14,fillstyle=solid,dotstyle=o](6.1250486,0.97557646)
\psdots[dotsize=0.14,fillstyle=solid,dotstyle=o](0.3050486,0.99557644)
\psdots[dotsize=0.14,fillstyle=solid,dotstyle=o](6.1250486,-1.0444236)
\psdots[dotsize=0.14,fillstyle=solid,dotstyle=o](0.32504863,-1.0444236)
\psline[linewidth=0.02cm](1.7809376,-1.0003124)(2.6409376,-1.0003124)
\psline[linewidth=0.02cm](2.7809374,-1.0203125)(3.6409376,-1.0203125)
\psline[linewidth=0.02cm](3.8009374,-1.0203125)(4.6409373,-1.0203125)
\psline[linewidth=0.02cm](1.8009375,0.9996875)(2.6809375,0.9996875)
\psline[linewidth=0.02cm](2.8209374,0.9996875)(3.6609375,0.9796875)
\psline[linewidth=0.02cm](3.7809374,0.9796875)(4.6409373,0.9996875)
\psline[linewidth=0.02cm,linestyle=dotted,dotsep=0.16cm](4.8009377,0.9996875)(6.0409374,0.9796875)
\psline[linewidth=0.02cm,linestyle=dotted,dotsep=0.16cm](4.7809377,-1.0003124)(6.0409374,-1.0203125)
\psline[linewidth=0.02cm,linestyle=dotted,dotsep=0.16cm](0.3609375,-1.0003124)(1.6609375,-1.0003124)
\psline[linewidth=0.02cm,linestyle=dotted,dotsep=0.16cm](0.3609375,1.0396875)(1.6809375,0.9996875)
\psline[linewidth=0.02cm](3.7609375,0.9196875)(6.0809374,-1.0003124)
\psline[linewidth=0.02cm](2.6809375,0.9596875)(0.3609375,-0.9803125)
\psline[linewidth=0.02cm](2.7609375,-0.9603125)(4.5409374,0.0196875)
\psline[linewidth=0.02cm](4.9009376,0.2196875)(6.0809374,0.9196875)
\psline[linewidth=0.02cm](3.7009375,-0.9603125)(3.3609376,-0.7403125)
\psline[linewidth=0.02cm](0.3409375,0.9596875)(1.6209375,0.2396875)
\psline[linewidth=0.02cm](1.8809375,0.0996875)(3.0809374,-0.6003125)
\usefont{T1}{ptm}{m}{n}
\rput(4.7323437,1.3096875){$-2$}
\usefont{T1}{ptm}{m}{n}
\rput(6.052344,1.2896875){$-2$}
\usefont{T1}{ptm}{m}{n}
\rput(4.7323437,-1.3303125){$-2$}
\usefont{T1}{ptm}{m}{n}
\rput(6.092344,-1.3303125){$-2$}
\usefont{T1}{ptm}{m}{n}
\rput(0.31234375,1.2896875){$-2$}
\usefont{T1}{ptm}{m}{n}
\rput(1.7123437,1.2696875){$-2$}
\usefont{T1}{ptm}{m}{n}
\rput(0.29234374,-1.3303125){$-2$}
\usefont{T1}{ptm}{m}{n}
\rput(1.6723437,-1.3103125){$-2$}
\usefont{T1}{ptm}{m}{n}
\rput(2.9023438,1.2896875){$-a_3$}
\usefont{T1}{ptm}{m}{n}
\rput(2.7223437,-1.2903125){$-a_2$}
\usefont{T1}{ptm}{m}{n}
\rput(3.7423437,1.2896875){$-a_1$}
\usefont{T1}{ptm}{m}{n}
\rput(3.7423437,-1.2903125){$-a_4$}
\psline[linewidth=0.02cm](4.7009373,1.4996876)(4.7009373,1.7196875)
\psline[linewidth=0.02cm](4.7009373,1.7196875)(6.1209373,1.7196875)
\psline[linewidth=0.02cm](6.1009374,1.7196875)(6.1009374,1.4996876)
\psline[linewidth=0.02cm](0.3009375,1.4796875)(0.3009375,1.6996875)
\psline[linewidth=0.02cm](0.3009375,1.6996875)(1.7209375,1.6996875)
\psline[linewidth=0.02cm](1.7209375,1.6996875)(1.7209375,1.4796875)
\psline[linewidth=0.02cm](0.3209375,-1.7403125)(0.3209375,-1.5203125)
\psline[linewidth=0.02cm](0.3209375,-1.7203125)(1.7409375,-1.7203125)
\psline[linewidth=0.02cm](1.7209375,-1.5203125)(1.7209375,-1.7403125)
\psline[linewidth=0.02cm](4.7009373,-1.7403125)(4.7009373,-1.5203125)
\psline[linewidth=0.02cm](4.7009373,-1.7203125)(6.1209373,-1.7203125)
\psline[linewidth=0.02cm](6.1009374,-1.5203125)(6.1009374,-1.7403125)
\usefont{T1}{ptm}{m}{n}
\rput(5.3923435,1.9896874){$a_2-1$}
\usefont{T1}{ptm}{m}{n}
\rput(5.3723435,-1.9703125){$a_1-1$}
\usefont{T1}{ptm}{m}{n}
\rput(1.0323437,-1.9303125){$a_3-1$}
\usefont{T1}{ptm}{m}{n}
\rput(0.9923437,2.0096874){$a_4-1$}
\end{pspicture}
}

\bigskip\noindent Here, $\bullet$  means a $(-1)$-curve. By contracting the
two maximal linear chains of rational curves denoted by white
vertices, $$Z(a_1, a_2, a_3, a_4)\to T(a_1, a_2, a_3, a_4),$$ we
get a rational $\mathbb{Q}$-homology projective plane $T(a_1, a_2,
a_3, a_4)$ with two cyclic quotient singularities. It has the same
singularities as Koll\'ar's example $X(a_1, a_2, a_3, a_4)$ when
$w^*=1$.

\begin{proposition}
\begin{enumerate}
    \item If $a_1, a_2, a_3, a_4\ge 3$ and $a_i >3$ for some
        $i$, then $K_T$ is ample.
    \item If $a_1=a_2=a_3=a_4= 3$, then $K_T$ is
        numerically trivial.
    \end{enumerate}
\end{proposition}

\begin{proof}
Let $E:=E_1$ be the $(-1)$-curve meeting the component of
self-intersection $-a_1$ of the upper chain and the rightmost
component of the bottom chain. By
 Lemma \ref{ample}, we see that

$$\begin{array}{l}
 E f^*(K_T) \\
   = 1 - \dfrac{|[\overset{a_2-1}{\overbrace{2,\ldots, 2}}]|
 + |[\overset{a_4-1}{\overbrace{2,\ldots, 2}}, a_3]|}{|[2*(a_4-1),a_3,a_1,2*(a_2-1)]|}
 - \dfrac{1 + |[\overset{a_3-1}{\overbrace{2,\ldots, 2}}, a_2, a_4, \overset{a_1-2}{\overbrace{2,\ldots, 2}}]|}{|[2*(a_3-1), a_2,a_4,2*(a_1-1)]|}.
  \end{array}$$
By a direct computation,
$$|[\overset{a_2-1}{\overbrace{2,\ldots, 2}}]|
 + |[\overset{a_4-1}{\overbrace{2,\ldots, 2}},
 a_3]|=a_2+(a_3a_4-a_4+1).$$
By Lemma \ref{cf},
 $$\begin{array}{ll} &|[\overset{a_3-1}{\overbrace{2,\ldots, 2}}, a_2, a_4, \overset{a_1-2}{\overbrace{2,\ldots,
 2}}]|\\ &\quad= |[2*(a_3-1), a_2,a_4,2*(a_1-2)]|\\
 &\quad=(a_1-1)a_2a_3a_4-(a_1-1)a_2a_3-a_4(a_1-1)a_3+(a_1-1)a_4\\
 &\quad\qquad + a_2a_3-(a_1-1)-a_3+1\\
 &\quad =|[2*(a_3-1), a_2,a_4,2*(a_1-1)]|-(a_2a_3a_4-a_2a_3-a_3a_4+a_4-1).
\end{array}$$

\bigskip\noindent
Thus
$$\begin{array}{l}
 E f^*(K_T)  = -
 \dfrac{a_3a_4+a_2-a_4+1}{|[2*(a_4-1),a_3,a_1,2*(a_2-1)]|}
 +\dfrac{a_2a_3a_4-a_2a_3-a_3a_4+a_4-2}{|[2*(a_3-1), a_2,a_4,2*(a_1-1)]|}\\
\quad =
\dfrac{(a_2a_3a_4-a_3a_4+a_4-1)\{(a_1-1)(a_2-1)(a_3-1)(a_4-1) -
a_1a_3-a_2a_4+2\}}
 {|[2*(a_4-1),a_3,a_1,2*(a_2-1)]|\cdot |[2*(a_3-1),
 a_2,a_4,2*(a_1-1)]|}.
 \end{array}$$

\bigskip\noindent
 Note first that $$a_2a_3a_4-a_3a_4+a_4-1>0$$ if all $a_i \geq 2$.

Assume that $a_i \geq 3$ for every $i=1,2,3,4$. Then, it is easy
to see that

\bigskip
 \begin{itemize}
 \item $a_1a_3 \geq 2(a_1 + a_3) -3,$ where the equality holds iff
 $a_1=a_3=3$,
 \item $a_2a_4 \geq 2(a_2 + a_4) -3,$ where the equality holds iff
 $a_2=a_4=3$,
 \item $(a_1-1)(a_3-1) \cdot (a_2-1)(a_4-1) \geq 2\{(a_1-1)(a_3-1) + (a_2-1)(a_4-1)\}$, where the equality holds iff
 $a_1=a_2=a_3=a_4=3$.
 \end{itemize}

 \bigskip\noindent
Thus
 $$\begin{array}{l}
  (a_1-1)(a_2-1)(a_3-1)(a_4-1) - a_1a_3-a_2a_4+2\\
 \quad\geq 2\{(a_1-1)(a_3-1) + (a_2-1)(a_4-1)\} - a_1a_3-a_2a_4+2\\
 \quad= a_1 a_3 - 2(a_1+a_3)+a_2a_4-2(a_2+a_4)+6\\
 \quad\geq 0.
 \end{array}$$
Here, both inequalities become equalities iff $a_1=a_2=a_3=a_4=3$.
Now we apply Lemma \ref{ample2} to get the assertions.
\end{proof}
This completes the proof of Theorem \ref{two}.

\begin{remark} It is easy to check the following:
\begin{enumerate}
\item When $a_i=a_j=2$ for $\{i,j\}\in \{\{1,2\}, \{1,4\},
\{2, 3\},\{3,4\}\}$, $K_T$ is ample iff $a_k, a_l\ge 6$ or $a_k=5,
a_l\ge 7$ or $a_k=4, a_l\ge 10$, where $\{i,j, k,l\}=
\{1,2,3,4\}$.
\item When $a_i=a_j=2$ for $\{i,j\}\in \{\{1,3\}, \{2,4\}\}$, $-K_T$ is ample for all $a_k, a_l\ge
2$, where $\{i,j, k,l\}= \{1,2,3,4\}$.
\end{enumerate}
\end{remark}

\section{Examples with one cyclic singularity}

In this section, we construct a series of new rational
$\mathbb{Q}$-homology projective planes with $K_S$ ample starting
from  a different configuration of curves in $\mathbb{P}^2$.

Consider the following configuration of 4 lines and a nodal cubic
curve in $\mathbb{P}^2$.

\bigskip
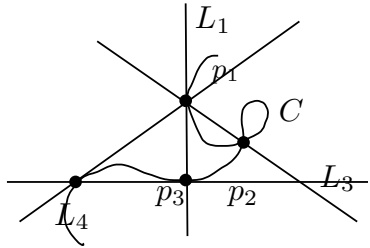
\scalebox{1.2} % Change this value to rescale the drawing.
{
\begin{pspicture}(-2.5, -1.5)(5,1.5)
\psline[linewidth=0.02cm](0.0,-0.64)(4.0,-0.64)
\psline[linewidth=0.02cm](1.98,1.34)(2.0,-1.24)
\psline[linewidth=0.02cm](0.14,-1.08)(3.24,1.14)
\psline[linewidth=0.02cm](3.88,-1.08)(1.0,0.92)
\psdots[dotsize=0.14](1.98,0.26) \psdots[dotsize=0.14](2.62,-0.2)
\psdots[dotsize=0.14](1.98,-0.62)
\psdots[dotsize=0.14](0.76,-0.64) \pscustom[linewidth=0.02] {
\newpath
\moveto(2.34,0.76) \lineto(2.29,0.76)
\curveto(2.265,0.76)(2.225,0.745)(2.21,0.73)
\curveto(2.195,0.715)(2.165,0.68)(2.15,0.66)
\curveto(2.135,0.64)(2.11,0.595)(2.1,0.57)
\curveto(2.09,0.545)(2.07,0.495)(2.06,0.47)
\curveto(2.05,0.445)(2.03,0.405)(2.02,0.39)
\curveto(2.01,0.375)(2.0,0.335)(2.0,0.31)
\curveto(2.0,0.285)(2.005,0.235)(2.01,0.21)
\curveto(2.015,0.185)(2.03,0.135)(2.04,0.11)
\curveto(2.05,0.085)(2.07,0.035)(2.08,0.01)
\curveto(2.09,-0.015)(2.11,-0.06)(2.12,-0.08)
\curveto(2.13,-0.1)(2.15,-0.14)(2.16,-0.16)
\curveto(2.17,-0.18)(2.205,-0.21)(2.23,-0.22)
\curveto(2.255,-0.23)(2.305,-0.24)(2.33,-0.24)
\curveto(2.355,-0.24)(2.405,-0.24)(2.43,-0.24)
\curveto(2.455,-0.24)(2.5,-0.235)(2.52,-0.23)
\curveto(2.54,-0.225)(2.585,-0.205)(2.61,-0.19)
\curveto(2.635,-0.175)(2.685,-0.15)(2.71,-0.14)
\curveto(2.735,-0.13)(2.78,-0.105)(2.8,-0.09)
\curveto(2.82,-0.075)(2.845,-0.04)(2.85,-0.02)
\curveto(2.855,0.0)(2.87,0.04)(2.88,0.06)
\curveto(2.89,0.08)(2.895,0.125)(2.89,0.15)
\curveto(2.885,0.175)(2.865,0.215)(2.85,0.23)
\curveto(2.835,0.245)(2.795,0.26)(2.77,0.26)
\curveto(2.745,0.26)(2.695,0.245)(2.67,0.23)
\curveto(2.645,0.215)(2.615,0.175)(2.61,0.15)
\curveto(2.605,0.125)(2.6,0.075)(2.6,0.05)
\curveto(2.6,0.025)(2.605,-0.02)(2.61,-0.04)
\curveto(2.615,-0.06)(2.62,-0.11)(2.62,-0.14)
\curveto(2.62,-0.17)(2.61,-0.225)(2.6,-0.25)
\curveto(2.59,-0.275)(2.57,-0.315)(2.56,-0.33)
\curveto(2.55,-0.345)(2.52,-0.38)(2.5,-0.4)
\curveto(2.48,-0.42)(2.44,-0.45)(2.42,-0.46)
\curveto(2.4,-0.47)(2.36,-0.495)(2.34,-0.51)
\curveto(2.32,-0.525)(2.275,-0.555)(2.25,-0.57)
\curveto(2.225,-0.585)(2.175,-0.605)(2.15,-0.61)
\curveto(2.125,-0.615)(2.075,-0.62)(2.05,-0.62)
\curveto(2.025,-0.62)(1.975,-0.62)(1.95,-0.62)
\curveto(1.925,-0.62)(1.875,-0.62)(1.85,-0.62)
\curveto(1.825,-0.62)(1.775,-0.615)(1.75,-0.61)
\curveto(1.725,-0.605)(1.675,-0.59)(1.65,-0.58)
\curveto(1.625,-0.57)(1.58,-0.555)(1.56,-0.55)
\curveto(1.54,-0.545)(1.5,-0.53)(1.48,-0.52)
\curveto(1.46,-0.51)(1.415,-0.49)(1.39,-0.48)
\curveto(1.365,-0.47)(1.315,-0.455)(1.29,-0.45)
\curveto(1.265,-0.445)(1.215,-0.445)(1.19,-0.45)
\curveto(1.165,-0.455)(1.12,-0.47)(1.1,-0.48)
\curveto(1.08,-0.49)(1.04,-0.51)(1.02,-0.52)
\curveto(1.0,-0.53)(0.955,-0.545)(0.93,-0.55)
\curveto(0.905,-0.555)(0.86,-0.575)(0.84,-0.59)
\curveto(0.82,-0.605)(0.785,-0.64)(0.77,-0.66)
\curveto(0.755,-0.68)(0.725,-0.72)(0.71,-0.74)
\curveto(0.695,-0.76)(0.675,-0.805)(0.67,-0.83)
\curveto(0.665,-0.855)(0.655,-0.905)(0.65,-0.93)
\curveto(0.645,-0.955)(0.64,-1.005)(0.64,-1.03)
\curveto(0.64,-1.055)(0.65,-1.105)(0.66,-1.13)
\curveto(0.67,-1.155)(0.7,-1.205)(0.72,-1.23)
\curveto(0.74,-1.255)(0.775,-1.295)(0.79,-1.31)
\curveto(0.805,-1.325)(0.83,-1.34)(0.84,-1.34)
\curveto(0.85,-1.34)(0.85,-1.33)(0.84,-1.32) }
\usefont{T1}{ptm}{m}{n}
\rput(3.6614063,-0.6){$L_3$}
\usefont{T1}{ptm}{m}{n}
\rput(3.0814063,1.96){$L_2$}
\usefont{T1}{ptm}{m}{n}
\rput(2.2814063,1.16){$L_1$}
\usefont{T1}{ptm}{m}{n}
\rput(3.1514063,0.16){$C$}
\usefont{T1}{ptm}{m}{n}
\rput(2.4275,0.55){\small $p_1$}
\usefont{T1}{ptm}{m}{n}
\rput(1.8214062,-0.8){$p_3$}
\usefont{T1}{ptm}{m}{n}
\rput(2.6214062,-0.8){$p_2$}
\usefont{T1}{ptm}{m}{n}
\rput(0.72140627,-1.02){$L_4$}
\end{pspicture}
}

%CLARIFY THIS PART%%%%%%%%%%%%%%%%%%%%%%%%%%%%%%%%%%%%%%%%%%
\noindent The existence of the configuration can be checked as
follows. Consider a plane nodal cubic curve $$C: y^2 = x^3 + x^2$$
on $\mathbb{C}^2$. Let $$L_1: y = ax$$ be a line passing through
the node of  $C$. If $a\neq \pm 1$, $L_1$ passes through $C$ at a
point $p_1$ different from the origin. For $i=2,3,4$, we
recursively define $L_i$ as the tangent line of $C$ at $p_{i-1}$,
and $p_i$ as the another intersection point of $L_i$ and $C$. For
a generic choice of $a$, we may assume that none of $p_1, p_2,
p_3$ is a flex of $C$, i.e., we may assume that $p_1\neq p_2,
p_2\neq p_3, p_3\neq p_4$, Now some calculation shows that there
is a suitable number $a\neq \pm 1$ such that $p_1 = p_4$, and
$p_1\neq p_2, p_2\neq p_3, p_3\neq p_1$. This shows the existence
of such a configuration.
%%%%%%%%%%%%%%%%%%%%%%%%%%%%%%%%%%%%%%%%%%

Now, by blowing up the node once and the three other marked points
three times each, we get a rational surface $Z(2)$ with the
following configuration of $15$ rational curves

\bigskip
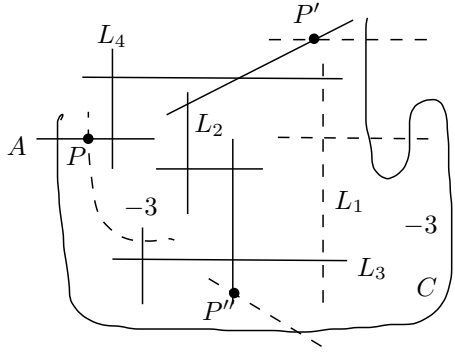
\scalebox{1} % Change this value to rescale the drawing.
{
\begin{pspicture}(-3, -2.3)(5,2.3)
\psline[linewidth=0.02cm](0.61,1.3407812)(4.07,1.3207812)
\psline[linewidth=0.02cm](2.01,1.1407813)(2.01,-0.47921875)
\psline[linewidth=0.02cm](1.59,0.12078125)(3.01,0.12078125)
\psline[linewidth=0.02cm](1.01,-1.0792187)(4.13,-1.0992187)
\psline[linewidth=0.02cm](1.01,1.7207812)(1.01,0.12078125)
\psline[linewidth=0.02cm](1.41,-0.6592187)(1.41,-1.6792188)
\psline[linewidth=0.02cm](0.0,0.5207813)(1.57,0.5207813)
\psline[linewidth=0.02cm](2.61,0.5207813)(2.61,-1.6792188)
\psline[linewidth=0.02cm](1.73,0.8407813)(4.21,2.1007812)
\psline[linewidth=0.02cm,linestyle=dashed,dash=0.16cm
0.16cm](3.81,1.5207813)(3.81,-1.6592188)
\psline[linewidth=0.02cm,linestyle=dashed,dash=0.16cm
0.16cm](3.09,1.8407812)(5.19,1.8407812)
\psline[linewidth=0.02cm,linestyle=dashed,dash=0.16cm
0.16cm](3.21,0.54078126)(5.21,0.5207813)
\psline[linewidth=0.02cm,linestyle=dashed,dash=0.16cm
0.16cm](2.29,-1.3392187)(3.79,-2.2392187)
\pscustom[linewidth=0.02,linestyle=dashed,dash=0.16cm 0.16cm] {
\newpath
\moveto(0.69,0.88078123) \lineto(0.69,0.60078126)
\curveto(0.69,0.46078125)(0.7,0.24578124)(0.71,0.17078125)
\curveto(0.72,0.09578125)(0.745,-0.09921875)(0.76,-0.21921875)
\curveto(0.775,-0.33921874)(0.825,-0.50421876)(0.86,-0.5492188)
\curveto(0.895,-0.59421873)(0.985,-0.68421876)(1.04,-0.7292187)
\curveto(1.095,-0.77421874)(1.225,-0.82421875)(1.3,-0.82921875)
\curveto(1.375,-0.83421874)(1.48,-0.83921874)(1.51,-0.83921874)
\curveto(1.54,-0.83921874)(1.59,-0.83421874)(1.61,-0.82921875)
\curveto(1.63,-0.82421875)(1.69,-0.81421876)(1.73,-0.80921876)
\curveto(1.77,-0.80421877)(1.82,-0.81421876)(1.83,-0.82921875) }
\usefont{T1}{ptm}{m}{n}
\rput(5.105625,-0.62921876){$-3$}
\usefont{T1}{ptm}{m}{n}
\rput(1.3856249,-0.38921875){$-3$}
\psdots[dotsize=0.14](3.69,1.8507813)
\usefont{T1}{ptm}{m}{n}
\rput(3.595625,2.2007812){$P'$}
\psdots[dotsize=0.14](2.63,-1.5307813)
\usefont{T1}{ptm}{m}{n}
\rput(2.43,-1.75){$P''$} \psdots[dotsize=0.14](0.69,0.5207813)
\pscustom[linewidth=0.02] {
\newpath
\moveto(4.4,2.1307812) \lineto(4.39,2.0907812)
\curveto(4.385,2.0707812)(4.38,2.0157812)(4.38,1.9807812)
\curveto(4.38,1.9457812)(4.38,1.8807813)(4.38,1.8507812)
\curveto(4.38,1.8207812)(4.38,1.7657813)(4.38,1.7407813)
\curveto(4.38,1.7157812)(4.38,1.6607813)(4.38,1.6307813)
\curveto(4.38,1.6007812)(4.38,1.5357813)(4.38,1.5007813)
\curveto(4.38,1.4657812)(4.38,1.4007813)(4.38,1.3707813)
\curveto(4.38,1.3407812)(4.38,1.2757813)(4.38,1.2407813)
\curveto(4.38,1.2057812)(4.38,1.1407813)(4.38,1.1107812)
\curveto(4.38,1.0807812)(4.38,1.0207813)(4.38,0.99078125)
\curveto(4.38,0.9607813)(4.38,0.90578127)(4.38,0.88078123)
\curveto(4.38,0.85578126)(4.38,0.79578125)(4.38,0.7607812)
\curveto(4.38,0.72578126)(4.39,0.65578127)(4.4,0.62078124)
\curveto(4.41,0.5857813)(4.425,0.5257813)(4.43,0.50078124)
\curveto(4.435,0.47578126)(4.445,0.42078125)(4.45,0.39078125)
\curveto(4.455,0.36078125)(4.46,0.29078126)(4.46,0.25078124)
\curveto(4.46,0.21078125)(4.475,0.14578125)(4.49,0.12078125)
\curveto(4.505,0.09578125)(4.545,0.05078125)(4.57,0.03078125)
\curveto(4.595,0.01078125)(4.645,-0.00921875)(4.67,-0.00921875)
\curveto(4.695,-0.00921875)(4.755,0.01578125)(4.79,0.04078125)
\curveto(4.825,0.06578125)(4.875,0.11578125)(4.89,0.14078125)
\curveto(4.905,0.16578124)(4.93,0.22578125)(4.94,0.26078126)
\curveto(4.95,0.29578125)(4.96,0.36078125)(4.96,0.39078125)
\curveto(4.96,0.42078125)(4.96,0.48578125)(4.96,0.5207813)
\curveto(4.96,0.55578125)(4.96,0.61578125)(4.96,0.6407812)
\curveto(4.96,0.66578126)(4.96,0.7157813)(4.96,0.74078125)
\curveto(4.96,0.7657812)(4.97,0.8257812)(4.98,0.86078125)
\curveto(4.99,0.8957813)(5.025,0.9507812)(5.05,0.97078127)
\curveto(5.075,0.99078125)(5.135,1.0207813)(5.17,1.0307813)
\curveto(5.205,1.0407813)(5.27,1.0457813)(5.3,1.0407813)
\curveto(5.33,1.0357813)(5.375,1.0057813)(5.39,0.98078126)
\curveto(5.405,0.9557812)(5.43,0.8957813)(5.44,0.86078125)
\curveto(5.45,0.8257812)(5.46,0.75578123)(5.46,0.72078127)
\curveto(5.46,0.68578124)(5.465,0.60578126)(5.47,0.56078124)
\curveto(5.475,0.5157812)(5.48,0.44578126)(5.48,0.42078125)
\curveto(5.48,0.39578125)(5.48,0.33078125)(5.48,0.29078126)
\curveto(5.48,0.25078124)(5.48,0.18078125)(5.48,0.15078124)
\curveto(5.48,0.12078125)(5.485,0.05578125)(5.49,0.02078125)
\curveto(5.495,-0.01421875)(5.5,-0.08421875)(5.5,-0.11921875)
\curveto(5.5,-0.15421875)(5.505,-0.22421876)(5.51,-0.25921875)
\curveto(5.515,-0.29421875)(5.52,-0.35921875)(5.52,-0.38921875)
\curveto(5.52,-0.41921875)(5.52,-0.47421876)(5.52,-0.49921876)
\curveto(5.52,-0.52421874)(5.52,-0.58921874)(5.52,-0.62921876)
\curveto(5.52,-0.6692188)(5.52,-0.7392188)(5.52,-0.76921874)
\curveto(5.52,-0.7992188)(5.52,-0.86921877)(5.52,-0.9092187)
\curveto(5.52,-0.94921875)(5.52,-1.0292188)(5.52,-1.0692188)
\curveto(5.52,-1.1092187)(5.52,-1.1842188)(5.52,-1.2192187)
\curveto(5.52,-1.2542187)(5.515,-1.3142188)(5.51,-1.3392187)
\curveto(5.505,-1.3642187)(5.49,-1.4292188)(5.48,-1.4692187)
\curveto(5.47,-1.5092187)(5.445,-1.5792187)(5.43,-1.6092187)
\curveto(5.415,-1.6392188)(5.375,-1.6992188)(5.35,-1.7292187)
\curveto(5.325,-1.7592187)(5.27,-1.8092188)(5.24,-1.8292187)
\curveto(5.21,-1.8492187)(5.15,-1.8842187)(5.12,-1.8992188)
\curveto(5.09,-1.9142188)(5.035,-1.9342188)(5.01,-1.9392188)
\curveto(4.985,-1.9442188)(4.915,-1.9492188)(4.87,-1.9492188)
\curveto(4.825,-1.9492188)(4.75,-1.9492188)(4.72,-1.9492188)
\curveto(4.69,-1.9492188)(4.615,-1.9592187)(4.57,-1.9692187)
\curveto(4.525,-1.9792187)(4.445,-1.9892187)(4.41,-1.9892187)
\curveto(4.375,-1.9892187)(4.31,-1.9892187)(4.28,-1.9892187)
\curveto(4.25,-1.9892187)(4.195,-1.9942187)(4.17,-1.9992187)
\curveto(4.145,-2.0042188)(4.09,-2.0092187)(4.06,-2.0092187)
\curveto(4.03,-2.0092187)(3.975,-2.0092187)(3.95,-2.0092187)
\curveto(3.925,-2.0092187)(3.87,-2.0142188)(3.84,-2.0192187)
\curveto(3.81,-2.0242188)(3.745,-2.0342188)(3.71,-2.0392187)
\curveto(3.675,-2.0442188)(3.605,-2.0492187)(3.57,-2.0492187)
\curveto(3.535,-2.0492187)(3.465,-2.0492187)(3.43,-2.0492187)
\curveto(3.395,-2.0492187)(3.325,-2.0442188)(3.29,-2.0392187)
\curveto(3.255,-2.0342188)(3.18,-2.0292187)(3.14,-2.0292187)
\curveto(3.1,-2.0292187)(3.015,-2.0292187)(2.97,-2.0292187)
\curveto(2.925,-2.0292187)(2.845,-2.0292187)(2.81,-2.0292187)
\curveto(2.775,-2.0292187)(2.705,-2.0292187)(2.67,-2.0292187)
\curveto(2.635,-2.0292187)(2.57,-2.0292187)(2.54,-2.0292187)
\curveto(2.51,-2.0292187)(2.45,-2.0242188)(2.42,-2.0192187)
\curveto(2.39,-2.0142188)(2.33,-2.0092187)(2.3,-2.0092187)
\curveto(2.27,-2.0092187)(2.21,-2.0092187)(2.18,-2.0092187)
\curveto(2.15,-2.0092187)(2.09,-2.0092187)(2.06,-2.0092187)
\curveto(2.03,-2.0092187)(1.96,-2.0092187)(1.92,-2.0092187)
\curveto(1.88,-2.0092187)(1.8,-2.0092187)(1.76,-2.0092187)
\curveto(1.72,-2.0092187)(1.65,-2.0092187)(1.62,-2.0092187)
\curveto(1.59,-2.0092187)(1.53,-2.0092187)(1.5,-2.0092187)
\curveto(1.47,-2.0092187)(1.395,-2.0092187)(1.35,-2.0092187)
\curveto(1.305,-2.0092187)(1.225,-1.9992187)(1.19,-1.9892187)
\curveto(1.155,-1.9792187)(1.085,-1.9692187)(1.05,-1.9692187)
\curveto(1.015,-1.9692187)(0.945,-1.9642187)(0.91,-1.9592187)
\curveto(0.875,-1.9542187)(0.805,-1.9342188)(0.77,-1.9192188)
\curveto(0.735,-1.9042188)(0.67,-1.8692187)(0.64,-1.8492187)
\curveto(0.61,-1.8292187)(0.56,-1.7842188)(0.54,-1.7592187)
\curveto(0.52,-1.7342187)(0.485,-1.6842188)(0.47,-1.6592188)
\curveto(0.455,-1.6342187)(0.44,-1.5742188)(0.44,-1.5392188)
\curveto(0.44,-1.5042187)(0.435,-1.4242188)(0.43,-1.3792187)
\curveto(0.425,-1.3342187)(0.415,-1.2542187)(0.41,-1.2192187)
\curveto(0.405,-1.1842188)(0.395,-1.1192187)(0.39,-1.0892187)
\curveto(0.385,-1.0592188)(0.375,-0.99921876)(0.37,-0.96921873)
\curveto(0.365,-0.93921876)(0.36,-0.87421876)(0.36,-0.83921874)
\curveto(0.36,-0.80421877)(0.355,-0.7342188)(0.35,-0.69921875)
\curveto(0.345,-0.6642187)(0.335,-0.59421873)(0.33,-0.55921876)
\curveto(0.325,-0.52421874)(0.315,-0.44921875)(0.31,-0.40921876)
\curveto(0.305,-0.36921874)(0.3,-0.29921874)(0.3,-0.26921874)
\curveto(0.3,-0.23921876)(0.295,-0.16921875)(0.29,-0.12921876)
\curveto(0.285,-0.08921875)(0.28,-0.01921875)(0.28,0.01078125)
\curveto(0.28,0.04078125)(0.28,0.10578125)(0.28,0.14078125)
\curveto(0.28,0.17578125)(0.28,0.24078125)(0.28,0.27078125)
\curveto(0.28,0.30078125)(0.28,0.36078125)(0.28,0.39078125)
\curveto(0.28,0.42078125)(0.28,0.48078126)(0.28,0.5107812)
\curveto(0.28,0.54078126)(0.28,0.60078126)(0.28,0.63078123)
\curveto(0.28,0.66078126)(0.28,0.7157813)(0.28,0.74078125)
\curveto(0.28,0.7657812)(0.3,0.81078124)(0.32,0.8307812)
\curveto(0.34,0.85078126)(0.355,0.85078126)(0.35,0.8307812)
\curveto(0.345,0.81078124)(0.33,0.78078127)(0.32,0.7707813) }
\usefont{T1}{ptm}{m}{n}
\rput(0.525625,0.3007812){$P$}
\usefont{T1}{ptm}{m}{n}
\rput(1,1.8907812){$L_4$}
\usefont{T1}{ptm}{m}{n}
\rput(5.1914063,-1.4492188){$C$}
\usefont{T1}{ptm}{m}{n}
\rput(2.3,0.7107813){$L_2$}
\usefont{T1}{ptm}{m}{n}
\rput(4.461406,-1.2){$L_3$}
\usefont{T1}{ptm}{m}{n}
\rput(4.1614063,-0.30921876){$L_1$}
\usefont{T1}{ptm}{m}{n}
\rput(-0.26234376,0.43078123){$A$}
\end{pspicture}
}

\bigskip\noindent
Here, $C$ and $L_i$ are the proper transforms of $C$ and $L_i$, a
dotted curve is a $(-1)$-curve and a solid curve is a $(-2)$-curve
if it is not specified as a $(-3)$-curve. The surface $Z(2)$
contains the following Hirzebruch-Jung string of rational curves

$$\underset{C}{\overset{-3}\circ}-\underset{A}{\overset{-2}\circ}-\underset{L_4}{\overset{-2}\circ}
-\underset{}{\overset{-2}\circ}-\underset{}{\overset{-2}\circ}-\underset{L_2}{\overset{-2}\circ}
-\underset{}{\overset{-2}\circ}-\underset{}{\overset{-2}\circ}-\underset{L_3}{\overset{-2}\circ}
-\underset{}{\overset{-3}\circ}.$$

\noindent  Blowing up $(b-2)$ times the marked point $P$, we get a
surface $Z(b)$ with the following Hirzebruch-Jung string of
rational curves

$$\underset{C}{\overset{-3}\circ}-\underset{A}{\overset{-b}\circ}-\underset{L_4}{\overset{-2}\circ}
-\underset{}{\overset{-2}\circ}-\underset{}{\overset{-2}\circ}-\underset{L_2}{\overset{-2}\circ}
-\underset{}{\overset{-2}\circ}-\underset{}{\overset{-2}\circ}-\underset{L_3}{\overset{-2}\circ}
-\underset{}{\overset{-3}\circ}-\overset{b-2}{\overbrace{\underset{}{\overset{-2}\circ}-\cdots
-\underset{B}{\overset{-2}\circ}}}.$$

\noindent Now by contracting these rational curves,
$$Z(b)\to S(b),$$
we get a rational $\mathbb{Q}$-homology projective plane $S(b)$
with a unique cyclic singularity of type
$$\frac{1}{27b^2-36b+4}(1, 9b^2-9b+1).$$

Let $E$ be the exceptional curve of $Z(b)\to Z(b-1)$, i.e., the
$(-1)$-curve with $E.A = E.B = 1$. Then by Lemma \ref{ample},
$$\begin{array}{lll}
 E f^*(K_S) &=& 1 -\dfrac{3 +  |[\overset{7}{\overbrace{2, \ldots, 2}},3, \overset{b-2}{\overbrace{2,\ldots, 2}}]|}
{|[3,b,\underset{7}{\underbrace{2, \ldots, 2}},3,\underset{b-2}{\underbrace{2, \ldots, 2}}]|}
 - \dfrac{1 + |[3,b,\overset{7}{\overbrace{2, \ldots, 2}},3, \overset{b-3}{\overbrace{2,\ldots, 2}}]|}
 {|[3,b,\underset{7}{\underbrace{2, \ldots, 2}},3,
\underset{b-2}{\underbrace{2, \ldots, 2}}]|}\\
& = & 1 - \dfrac{3 + (9b-1)}{27b^2-36b+4} - \dfrac{1 + (27b^2-63b+37)}{27b^2-36b+4}\\
& = & \dfrac{18(b-2)}{27b^2-36b+4}.
\end{array}$$
Now Lemma \ref{ample2} completes the proof of Theorem \ref{one}.

\begin{remark} One can get more examples by blowing up not only the marked point
$P$ but also the marked point $P'$ or $P''$.
\begin{enumerate}
\item Blowing up $(b-2)$ times the marked point $P$ and $(c-2)$ times the marked point $P'$, we get a
surface $Z(b, c)$ with the following Hirzebruch-Jung string of
rational curves
$$[\overset{c-2}{\overbrace{2, \ldots, 2}},3, b, 2,2,c,2,2,2,2,3, \overset{b-2}{\overbrace{2,\ldots,
2}}].$$ Here $b, c\ge 2$. The resulting rational
$\mathbb{Q}$-homology projective plane $S(b,c)$ has ample
canonical class if $b$ and $c$ are not small.
\item Blowing up $(b-2)$ times the marked point $P$ and $(c-2)$ times the marked point $P''$, we get a
surface $W(b, c)$ with the following Hirzebruch-Jung string of
rational curves
$$[\overset{c-2}{\overbrace{2, \ldots, 2}},3, b, 2,2,2,2,2,c,2,3, \overset{b-2}{\overbrace{2,\ldots,
2}}].$$ Here $b, c\ge 2$. The resulting rational
$\mathbb{Q}$-homology projective plane $T(b,c)$ has ample
canonical class if $b$ and $c$ are not small.
\end{enumerate}
\end{remark}

\section{Examples with three cyclic singularities}

Consider the following configuration of 3 concurrent lines and a
conic on $\mathbb{P}^2$.

\bigskip

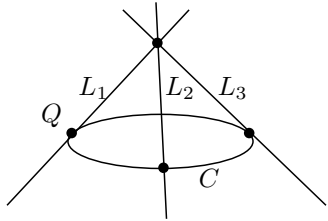
\scalebox{1} % Change this value to rescale the drawing.
{
\begin{pspicture}(-3.5, -1.5)(5,1.5)
\psellipse[linewidth=0.02,dimen=outer](2.04,-0.41)(1.24,0.37)
\psline[linewidth=0.02cm](0.0,-1.26)(2.42,1.34)
\psline[linewidth=0.02cm](4.24,-1.26)(1.54,1.34)
\psline[linewidth=0.02cm](1.98,1.44)(2.12,-1.44)
\psdots[dotsize=0.14](2.0,0.9) \psdots[dotsize=0.14](0.86,-0.3)
\psdots[dotsize=0.14](3.22,-0.3) \psdots[dotsize=0.14](2.08,-0.76)
\usefont{T1}{ptm}{m}{n}
\rput(0.5914062,-0.06){$Q$}
\usefont{T1}{ptm}{m}{n}
\rput(1.15,0.3){$L_1$}
\usefont{T1}{ptm}{m}{n}
\rput(2.3,0.3){$L_2$}
\usefont{T1}{ptm}{m}{n}
\rput(3.0014062,0.3){$L_3$}
\usefont{T1}{ptm}{m}{n}
\rput(2.6914062,-0.9){$C$}
\end{pspicture}
}

\bigskip\noindent
By blowing up the marked point $Q$ three times and the three other
marked points twice each, we get a rational surface $Z(2)$ with
the following configuration of $13$ rational curves.

\bigskip
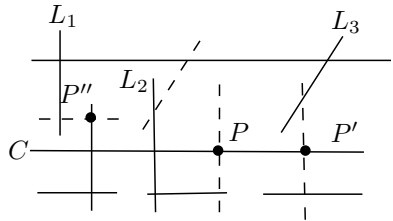
\scalebox{1} % Change this value to rescale the drawing.
{
\begin{pspicture}(-3.5, -1.5)(5,1.5)
\psline[linewidth=0.02cm](0.4309375,0.6907815)(5.2109375,0.6907815)
\psline[linewidth=0.02cm](0.4109375,-0.50921863)(4.8509374,-0.52921855)
\psline[linewidth=0.02cm](4.5709376,1.0107814)(3.7509375,-0.26921862)
\psline[linewidth=0.02cm,linestyle=dashed,dash=0.16cm
0.16cm](4.0309377,0.38921878)(4.0709376,-1.4292188)
\psline[linewidth=0.02cm](4.8309374,-1.0892185)(3.5109375,-1.0892185)
\psline[linewidth=0.02cm,linestyle=dashed,dash=0.16cm
0.16cm](2.6709375,0.9507814)(1.9109375,-0.22921862)
\psline[linewidth=0.02cm](2.0509374,0.4507814)(2.0805774,-1.3292185)
\psline[linewidth=0.02cm,linestyle=dashed,dash=0.16cm
0.16cm](2.9309375,0.3707814)(2.9297545,-1.3292185)
\psline[linewidth=0.02cm](0.8109375,1.0907812)(0.8109375,-0.30921862)
\psline[linewidth=0.02cm](1.2309375,0.1107814)(1.2309375,-1.3092185)
\psline[linewidth=0.02cm,linestyle=dashed,dash=0.16cm
0.16cm](0.5309375,-0.0892186)(1.5909375,-0.0892186)
\psline[linewidth=0.02cm](1.9704452,-1.0892185)(3.0309374,-1.070781)
\psdots[dotsize=0.14](2.9109375,-0.50921863)

\usefont{T1}{ptm}{m}{n}
\rput(4.595625,-0.25843734){$P'$}
\psdots[dotsize=0.14](4.073,-0.50921863)
\usefont{T1}{ptm}{m}{n}
\rput(1.03,0.25){$P''$} \psdots[dotsize=0.14](1.22,-0.069)

\usefont{T1}{ptm}{m}{n}
\rput(3.1865625,-0.25843734){$P$}
\psline[linewidth=0.02cm](1.5695312,-1.0699997)(0.50953126,-1.0699997)
\usefont{T1}{ptm}{m}{n}
%\rput(1.1665625,0.98156255){$-3$}
\usefont{T1}{ptm}{m}{n}
\rput(1.7923437,0.42000002){$L_2$}
\usefont{T1}{ptm}{m}{n}
\rput(0.26234376,-0.49999997){$C$}
\usefont{T1}{ptm}{m}{n}
\rput(0.8523437,1.28){$L_1$}
 \rput(4.612344,1.2){$L_3$}
%\rput(5.186562,-0.4984373){$-3$}
\end{pspicture}
}

\bigskip\noindent
Here, $C$ and $L_i$ are the proper transforms of $C$ and $L_i$, a
dotted line is a $(-1)$-curve and a solid line is a $(-2)$-curve
except $L_1$ which is a $(-3)$-curve. Note that $Z(2)$ contains
the following three Hirzebruch-Jung strings of rational curves

$$\overset{-2}\circ,\quad\quad
\underset{L_1}{\overset{-3}\circ}-\overset{-2}\circ-\underset{L_3}{\overset{-2}\circ},\quad\quad
\underset{}{\overset{-2}\circ}-\underset{}{\overset{-2}\circ}-\underset{C}{\overset{-2}\circ}
-\underset{L_2}{\overset{-2}\circ}-\underset{}{\overset{-2}\circ}.$$

\noindent Blowing up $(b-2)$ times the marked point $P$, we get a
rational surface $Z(b)$ with the following three Hirzebruch-Jung
strings of rational curves

$$\overset{-2}\circ,\quad\quad
\underset{L_1}{\overset{-3}\circ}-\overset{-2}\circ-\underset{L_3}{\overset{-2}\circ},\quad\quad
\underset{}{\overset{-2}\circ}-\underset{}{\overset{-2}\circ}-\underset{C}{\overset{-b}\circ}
-\overset{b}{\overbrace{\underset{L_2}{\overset{-2}\circ}-\cdots
-\underset{B}{\overset{-2}\circ}}}.$$

\noindent Contracting these  rational curves,
 $$Z(b)\to S(b),$$
 we get a
rational $\mathbb{Q}$-homology projective plane $S(b)$ with three
cyclic singularities of type
$$A_1, \,\,\, \frac{1}{7}(1,3),\,\,\, \frac{1}{3b^2-2b-2}(1, 2b^2-b-1).$$

Let $E$ be the exceptional curve of $Z(b)\to Z(b-1)$, i.e., the
$(-1)$-curve with $E.C = E.B = 1$. Then by Lemma \ref{ample},
$$ E f^*(K_S)  =  1 -\frac{(b+1)+3}{3b^2-2b-2}  -\frac{1 + (3b^2-5b+3)}{3b^2-2b-2}
 =   \frac{2(b-5)}{3b^2-2b-2}. $$
Now Lemma \ref{ample2} completes the proof of Theorem \ref{three}.

\begin{remark} One can get more examples by blowing up not only the marked point
$P$ but also the marked points $P'$ and $P''$.
\begin{enumerate}
\item Blowing up $(b-2)$ times the marked point $P$ and $c$ times the marked point $P''$, we get a
surface with the following three Hirzebruch-Jung strings of
rational curves
$$[2], \quad \quad [\overset{c}{\overbrace{2,\ldots,
2}},3,2,2], \quad\quad [2,2+c, b, \overset{b}{\overbrace{2,\ldots,
2}}].$$ Here $b\ge 2, c\ge 0$. The resulting rational
$\mathbb{Q}$-homology projective plane $V(b,c)$ has ample
canonical class if $b$ and $c$ are not small.
\item Blowing up $(b-2)$ times the marked point $P$,  once the marked point $P'$ and $c$ times the marked point $P''$, we get a
surface with the following two Hirzebruch-Jung strings of rational
curves
$$[\overset{c}{\overbrace{2,\ldots,
2}},3,2,2,2,2], \quad\quad [2,2+c, b+1,
\overset{b}{\overbrace{2,\ldots, 2}}].$$ Here $b\ge 2, c\ge 0$.
The resulting rational $\mathbb{Q}$-homology projective plane
$Y(b,c)$ has 2 cyclic singularities, and its canonical class is
ample if $b$ and $c$ are not small.
\end{enumerate}
\end{remark}
\medskip

%%%%%%%%%%%%%%%%%%%%%%%%%%%%%%%%%%%%%%%%%%%%%%%%%%%%%%%%%%%%%%%%%%%%%%%%
%\bibliographystyle{amsplain}

\end{document}